# A strongly polynomial-time algorithm for the general linear programming problem

Samuel Awoniyi
Department of Industrial and Manufacturing Engineering
FAMU-FSU College of Engineering
2525 Pottsdamer Street
Tallahassee, FL 32310
E-mail: awoniyi@eng.famu.fsu.edu
ORCID# 0000 0001 7102 6257

.

# Abstract

This article presents a strongly polynomial-time algorithm for the general linear programming problem. This algorithm is an implicit reduction procedure that works as follows. Primal and dual problems are combined into a special system of linear equations constrained by complementarity relations and non-negative variables. Each iteration of the algorithm consists of applying a pair of complementary Gauss-Jordan pivoting operations, guided by a necessary-condition lemma. The algorithm requires no more than $2(k + n)$ iterations, as there are only $k + n$ complementary pairs of columns to compare one-pair-at-a-time, where $k$ is the number of constraints and $n$ is the number of variables of given general linear programming problem. Numerical illustration is given that includes an instance of a classical problem of Klee and Minty and a problem of Beale.

# 1. Introduction

To find a strongly polynomial-time algorithm for the general linear programming (LP) problem is still an open problem [5,6,7,16]. Classical references on this topic include [1,2,8,9,10,12,14,15,16].

Well-known efficacy of Dantzig's simplex algorithm variants on real-world LP problems suggests that someone should eventually find a strongly polynomial-time algorithm for the general LP problem before fast computing makes it a moot point. Such an algorithm might not only explain the success of current simplex algorithm variants when applied to practical problems, but might also indicate ways to enhance existing commercial software packages that utilize simplex algorithm variants and some of the newer algorithms inspired by [9, 10].

This article presents a strongly polynomial-time algorithm for solving the general LP problem. This algorithm begins by utilizing basic LP duality theory to translate solving the general LP problem, having $k$ inequality constraints and $n$ variables, into solving a special system of equations in $R^{2(k+n)}$. Each iteration consists of two special Gauss-Jordan reduction pivoting instances, in accordance with a necessary condition lemma (Lemma 6.1) proved in this article. The algorithm stops in at most $2(k + n)$ iterations.

The rest of this article is organized as follows. Section 2 gives our problem statement; there, the general LP problem is translated into a constrained equation-solving problem denoted as (Eq). Section 3 states our strategy for solving problem (Eq), thereby indicating a point of departure relative to existing algorithms for solving the general LP problem. Step-by-step details of our algorithm are stated in Section 4, followed by a numerical illustration in Section 5. Section 6 proves Lemma 6.1 which is a necessary-condition lemma, relative to computing a solution for (Eq). A sufficient-condition corollary of Lemma 6.1 is also proved in Section 6. The algorithm's computational complexity is explained with two Claims and Lemma 7.1 in Section 7.

## 2. Problem statement

An equation solving problem, denoted as (Eq), is introduced here. Problem (Eq) is a primal-dual translation of the LP problem.

As notation in this article, vectors are column vectors unless otherwise indicated. Vectors will be denoted by lower-case letters, and matrices by upper-case letters. Superscript $T$ will denote vector or matrix transpose as usual, and $I_{(.)}$ is reserved for identity matrix of dimension indicated by $(.)$.

.

We assume the general LP problem to be given in Neumann symmetric form, (P) below:

.

$$\left\{ \begin{array}{ll} \text{maximize} & f^T x \\ \text{subject to:} & Ax \le b \\ & x \ge 0 \end{array} \right\} \cdots\cdots (P)$$

.

where $f$ is $n$-vector, $A$ is $k$-by-$n$ (numerical) matrix, $b$ is $k$-vector, and $x$ is $n$-vector of problem's variables.

From basic LP duality theory, solving (P) is equivalent to computing a $2(k+n)$-vector $z$ that solves the constrained system of linear equations (Eq) stated below:

.

$$\left\{ \begin{array}{l} Mz = q, \\ z_j z_{(k+n+j)} = 0, \ \text{for}\, j = 1, \ldots, k+n \\ z \ge 0 \end{array} \right\} \cdots\cdots (Eq)$$

.

where

$$M = \begin{pmatrix} O & A & I_{(k)} & O \\ -A^T & O & O & I_{(n)} \\ -b^T & f^T & o^T & o^T \end{pmatrix} \text{and } q = \begin{pmatrix} b \\ -f \\ o \end{pmatrix}$$

.

As a numerical illustration of this problem statement, suppose

.

$$f = \begin{pmatrix} -1 \\ 1 \end{pmatrix}, \ A = \begin{pmatrix} 1 & 1 \\ -1 & 0 \end{pmatrix}, \ b = \begin{pmatrix} 10 \\ -5 \end{pmatrix}$$

.

In this instance, we have

.

$$
M = \begin{bmatrix} 0 & 0 & 1 & 1 & 1 & 0 & 0 & 0 \\ 0 & 0 & -1 & 0 & 0 & 1 & 0 & 0 \\ -1 & 1 & 0 & 0 & 0 & 0 & 1 & 0 \\ -1 & 0 & 0 & 0 & 0 & 0 & 0 & 1 \\ -10 & 5 & -1 & 1 & 0 & 0 & 0 & 0 \end{bmatrix} \text{ and } q = \begin{bmatrix} 10 \\ -5 \\ 1 \\ -1 \\ 0 \end{bmatrix}
$$
.

Problem (Eq) is an instance of primal-dual formulation of the general LP problem (see [4] for example). A widely studied primal-dual formulation known as the Linear Complementarity Problem (LCP) (see [3, 4, 11], for example) is what one gets from (Eq) by not including the last equation of $Mz = q$.

As it is well-known that the LCP includes the general LP problem in terms of solution existence, one might then surmise that the last equation of $Mz = q$ is redundant information. But our Lemma 6.1 in this article shows that the last equation of $Mz = q$ is indeed indispensable for our algorithm to work as described in this article.

# 3. Informal statement of our strategy for solving (Eq)

We indicate here how the strategy of our algorithm differs from those of existing algorithms for solving the general LP problem.

One may regard any algorithm for solving the general LP problem as either an improvement-direction-following algorithm (henceforth abbreviated as $\underset{\rightarrow}{id}$ algorithm) or as a non-$\underset{\rightarrow}{id}$ algorithm. In each step, an $\underset{\rightarrow}{id}$ algorithm follows a direction that improves upon current approximate solution until some convergence is attained, whereas a non-$\underset{\rightarrow}{id}$ algorithm instead utilizes only linear algebra and combinatorial facts to "leap onto" a solution after a while, usually on account of avoiding a contradiction, or after a well-defined set is exhausted.

With that informal definition of $\underset{\rightarrow}{id}$ algorithms, one can see that existing variants of Dantzig's simplex algorithm, Khachian's algorithm and Karmarkar's algorithm are $\underset{\rightarrow}{id}$ algorithms (see [6, 7] for comprehensive references on those algorithms), whereas the Fourier-Motzkin Elimination (FME) method (see pp. 43-51 of [5] for an introduction to the FME method) is a non-$\underset{\rightarrow}{id}$ algorithm.

The most efficacious among existing algorithms for solving the general LP problem are polynomial-time algorithms, with worst-case running time depending also on data coding size, say $H$. That is only stating the well-known fact that none of those algorithms is a strongly polynomial-time algorithm without major data restrictions (see [2, 12, 15] for example).

Reasoning from our observation of published literature on solving LP problems, along with our own recent computational experience with iterative methods for solving linear equations [17], we recently surmised that the worst-case time of $\underset{\rightarrow}{id}$ algorithms for the general LP problem will necessarily depend on problem data size $H$.

We therefore thought that a strategy of implicit reduction of solution space dimension (similar to Gaussian reduction) might not depend on $H$, thereby possibly resulting in a strongly polynomial-time algorithm for the general LP problem. That lofty supposition is a key motivation for our non-$\underset{\rightarrow}{id}$ algorithm presented in this article, and the following is an informal statement of another aspect of that motivation.

From basic linear programming theory, it is clear that if the system (Eq) has a solution, then it has one for which the basis matrix (a submatrix of $M$) is nonsingular. It follows then that solving problem (Eq) reduces to correctly selecting columns (of $M$) forming a solution basis matrix, without one having to devote any special attention to solution feasibility or nearness from a solution, since a correct basis matrix would ordinarily yield an optimal feasible solution.

Towards selecting such a correct solution basis matrix, the complementary slackness constraints $z_j z_{(k+n+j)} = 0,\ j = 1,\dots,k+n$ in (Eq) are helpful because they imply that selecting column $j$ precludes selecting column $k+n+j$, and vice versa, thereby defining an underlying set that can be exhausted after some $k+n$ column selections. The algorithm presented in this article may be regarded as a procedure that works in that vein.

# 4. Step-by-step description of our algorithm

We describe details of our algorithm step-by-step, in a manner that can be readily translated into a computer program. A validation of this algorithm will be given in Sections 6 and 7 where we explain the algorithm's computational complexity.

Our algorithm is a non-$\underrightarrow{id}$ algorithm for solving (Eq), with each iteration consisting of a pair of special Gauss-Jordan (abbreviated as GJ) pivoting guided by a necessary-condition lemma. The main steps of our algorithm are as follows.

.

Step 1 - `Initialize` the algorithm (as described below)

Step 2 - `Stop,` if stopping condition (described below) is met; otherwise go to Step 3

Step 3 - `Execute next iteration` (as described below), and thereafter go back to Step 2 above

## 4.1 Initialization

Initialization consists of setting up an initial "organizer" or "tableau" for the algorithm's iterations. As notation henceforth, we let $[M\,q]$ denote the augmented matrix combining matrix $M$ and column vector $q$ ($M$ and $q$ as introduced in problem (Eq) in Section 2).

.

`Initialization operation`

.

Add the $(k+n+1)$-th row of $[M\,q]$ to every other row of $[M\,q]$, in order to facilitate needed (complementary) Gauss-Jordan (GJ) pivoting (with pivoting positions) on diagonal elements of $[M\,q]$ later in the algorithm. Let us temporarily denote the resultant matrix (that is, $[M\,q]$ after that operation) by $[\overline{M}\,\overline{q}]$.

As numerical illustration of that initialization operation, let us consider again our illustration example of Section 2. We then have

.

$$[\overline{M}\,\overline{q}] = \begin{bmatrix} -10 & 5 & 0 & 2 & 1 & 0 & 0 & 0 & 10 \\ -10 & 5 & -2 & 1 & 0 & 1 & 0 & 0 & -5 \\ -11 & 6 & -1 & 1 & 0 & 0 & 1 & 0 & 1 \\ -11 & 5 & -1 & 1 & 0 & 0 & 0 & 1 & -1 \\ -10 & 5 & -1 & 1 & 0 & 0 & 0 & 0 & 0 \end{bmatrix}$$

.

From basic linear algebra, one can see that the system of equations and inequalities associated with $[\overline{M}\ \overline{q}]$ is equivalent to (Eq) in terms of solution existence. Accordingly, in the interest of notation efficacy, we will generally write $[M\ q]$ in place of $[\overline{M}\ \overline{q}]$ or any of its equivalent systems that we will obtain (in the algorithm) through elementary row operations on $[M\ q]$.

As a part of this initialization, define $\Pi$ as the set of indices of columns of $M$ that are *known* to be utilized (in a solution basis matrix) or *known* not to be utilized by some solution of (Eq), assuming that (Eq) has a solution. The index set $\Pi$ is initialized as the empty set $\varnothing$. Later, $\Pi$ will be updated continually, as described under "Executing next iteration" below.

In the remainder of this article, it may be advisable to keep in view the following general form of the matrix $[M\ q]$ :

| $m_{1,1}$ | $m_{1,2}$ | $\cdots$ | $m_{1,2(k+n)}$ | $q_1$ |
|---|---|---|---|---|
| $\vdots$ | $\vdots$ | $\ddots$ | $\vdots$ | $\vdots$ |
| $m_{k+n,1}$ | $m_{k+n,2}$ | $\cdots$ | $m_{k+n,2(k+n)}$ | $q_{k+n}$ |
| $m_{k+n+1,1}$ | $m_{k+n+1,2}$ | $\cdots$ | $m_{k+n+1,2(k+n)}$ | $q_{k+n+1}$ |

## 4.2 Stopping

There are two types of stopping - the case when a solution for (Eq) is found, and the case when there is evidence that (Eq) has no solutions.

.

`Case 1: A solution of (Eq) is found`

.

A solution of (Eq) is indicated in $[M\ q]$ by having $q \geq 0$ with $q_{k+n+1} = 0$.

`Case 2: There is evidence that (Eq) has no solutions`

.

A lack of solutions for (Eq) is indicated by having $q_{k+n+1} > 0$ with all other elements non-positive (that is, $\leq 0$) in row $k+n+1$ of $[M\ q]$, that possibly after we first multiply row $k+n+1$ by $-1$ to have $q_{k+n+1} > 0$. A lack of solution may also be indicated through Step 1 or Step 4 of Section 4.3 – Executing next iteration.

## 4.3 Executing next iteration

Towards describing "Executing next iteration", we define, in (i),..., (iv) below, several concepts that will feature frequently in the remainder of this article.

.

`Definitions`

.

(i) A GJ pivoting in column $j$ of $[M\ q]$ is called *the complementary GJ pivoting in column j* if the pivoting position in $[M\ q]$ is (j,j) for $j \leq k+n$, or (j-k-n,j) for $j > k+n$.

(ii) For $j \leq k+n$, column $j$ is the *complement column* for column $k+n+j$, and vice versa. Also, the index $j$ may be refered to as complement for the index $k+n+j$, and vice versa.

(iii) Each iteration of the algorithm consists of two complementary GJ pivoting instances – a Minor Pivoting (abbreviated as *MinorP*) instance, when $q_{k+n+1} = 0$, and a Major Pivoting (abbreviated as *MajorP*) instance, when $q_{k+n+1} > 0$ (or, equivalently, $q_{k+n+1} < 0$). MinorP pivoting and MajorP pivoting instances will be described shortly in this section.

(iv) We will sometimes refer to a column of $M$ in $[M\ q]$, say $M^{(s)}$, as a *maximal column* if $m_{k+n+1,s} \geq m_{k+n+1,j}$ for $j = 1,\ldots,k+n$; that is, its $(k+n+1)$-th component (its element in the last row of $M$) is not smaller than that of any other column of $M$.

.

An illustration of the definitions above

.

The two concepts defined in (i) and (ii) are illustrated by the following table.

| * | | | | * | | | | $q_1$ |
|---|---|---|---|---|---|---|---|---|
| | * | | | | * | | | $q_2$ |
| | | * | | | | * | | $q_3$ |
| | | | * | | | | * | $q_4$ |
| last row (row 5) of $M$ | | | | | | | | $q_5$ |

This table displays a frame of $[M\ q]$ for $k = n = 2$. An asterisk in the (i,j)-th position indicates the complementary pivoting position in column $j$, for $j = 1,...,8$. In each row, there are asterisks in two positions, and the columns of the two positions are complement columns for each other; that is, in general, if (i,j) is one (with j < k+n), then (i,k+n+j) is the other.

.

We next describe "MinorP pivoting instance" and "MajorP pivoting instance". For each type of pivoting instance, the current iteration (of the algorithm) utilizes a four-step procedure to "select" a column of $M$, and, thereafter, to perform the desired complementary GJ pivoting in the selected column. Details of the two procedures are as follows.

.

MinorP pivoting instance (Here $q_{k+n+1} = 0$ and $q_i < 0$ some $i \leq k+n$)

If necessary, multiply the last row of $[M\ q]$ (that is, row k+n+1) by –1 to ensure that each negative component of $q$ (in $[M\ q]$) corresponds to a positive component of the last row of $M$. Claim 4 in the Appendix of this article and the article arXiv:2410.19350 explain the feasibility of doing that.

.

*Step 1: Let $\mathcal{L}$ be the ordered list of column indices $j$ having $m_{k+n+1,j} > 0$, in descending order of $m_{k+n+1,j} > 0$. Let L denote $\mathcal{L}\backslash\Pi$. The items in L are to be "picked up" one-at-a-time for processing. If L = $\varnothing$, then current iteration declares that (Eq) has no solutions by virtue of an LP strict complementarity theorem, as $\mathcal{L} \neq \varnothing$; otherwise, there are two cases to consider here.

.. Case 1.1: L has exactly 1 item in it – In this case, letting L = $\{\hat{j}\}$, current iteration (i) updates $\Pi$ by putting $\hat{j}$ and its complement index into $\Pi$, and (ii) performs desired MinorP pivoting in $\hat{j}$-th column.

.. Case 1.2: L does not have exactly 1 item in it – In this case, with L having at least 2 elements in it, current iteration goes to Step 2.

*Step 2: Current iteration processes the next not-yet-processed item in list L. There are two cases to consider here.

.. Case 2.1: Some items are still available in list L to be processed – In this case, current iteration "picks up" the next item in L, labels the corresponding column as "current potential column selection", and then goes to Step 3.

..Case 2.2: There are no not-yet-processed items in L – In this case, current iteration goes to Step 4 where "finalizing operations" will be performed on the set L.

*Step 3: Current iteration processes the given "current potential column selection". There are two cases.

.. Case 3.1: The given "current potential column selection" is not the complement column for a previous MajorP "column selection" – In this case, the "current potential column selection" is labelled as "column selection" for this instance of MinorP, and current iteration thereafter performs desired MinorP pivoting in the column just labelled as "column selection".

.. Case 3.2: The given "current potential column selection" is the complement column for a previous MajorP "column selection" – In this case, the iteration goes back to Step 2, to process the next not-yet-processed item in list L.

*Step 4: First, do Task 1, Task 2 & Task 3 below, for each $\widehat{j} \in$ L separately, until *either* a solution of (Eq) is obtained, *or* the current iteration is ended with $\widehat{j}$ put in the index set $\Pi$.

Task 1: As a probe, imagine that MinorP pivoting is to be performed in column $\widehat{j}$, followed by pertinent MajorP pivoting, say in column $w$.

Task 2: If that MajorP pivoting in $w$ yields a solution of (Eq), then terminate the algorithm at this juncture. Otherwise, go to Task 3.

Task 3: In place of the pivoting sequence imagined in Task 1, consider doing a complementary GJ pivoting in column $w$ first, followed by doing a complementary pivoting in column $\widehat{j}$. If that reversal would indicate that $\widehat{j}$ should be put in the index set $\Pi$, then do those two complementary GJ pivotings (that is, do complementary GJ pivoting in column $w$ first, followed by doing it in column $\widehat{j}$), put $\widehat{j}$ and its complement in $\Pi$, and terminate the current iteration at this juncture.

But, if doing Tasks 1, 2, 3 for elements of L exhausts the set L without yielding a solution of (Eq) as in Task 2, and without putting $\widehat{j}$ in $\Pi$ and ending the current iteration as in Task 3, then the algorithm is to be terminated at this juncture, along with the declaration that (Eq) has no solutions. (Claim 7.2 in Section 7 validates that conclusion).

.

MajorP pivoting instance (Here $q_{k+n+1} > 0$ or $q_{k+n+1} < 0$)

If $q_{k+n+1} < 0$, then implicitly multiply row $k + n + 1$ of $[M\ q]$ by $-1$, to have $q_{k+n+1} > 0$.

.

*Step 1: Let $\mathcal{L}$ be the ordered list of column indices $j$ having $m_{k+n+1,j} > 0$, in descending order of $m_{k+n+1,j} > 0$. Let L denote $\mathcal{L}\backslash\Pi$. The items in L are to be "picked up" one-at-a-time for processing. If L $= \varnothing$, then current iteration declares that (Eq) has no solutions by virtue of an LP complementarity theorem, as $\mathcal{L} \neq \varnothing$; otherwise, there are two cases to consider here.

.. Case 1.1: L has exactly 1 item in it – In this case, letting L $= \{\widehat{j}\}$, current iteration (i) updates $\Pi$ by putting $\widehat{j}$ and its complement index into $\Pi$, and (ii) performs desired MajorP pivoting in $\widehat{j}$-th column.

.. Case 1.2: L does not have exactly 1 item in it – In this case, with L having at least 2 elements in it, current iteration goes to Step 2.

*Step 2: Current iteration processes the next not-yet-processed item in list L. There are two cases.

.. Case 2.1: Some items are still available in list L to be processed – In this case, current iteration "picks up" the next available item in L, labels the corresponding column as "current potential column selection" and then goes to Step 3.

.. Case 2.2: There are no not-yet-processed items in L – In this case, current iteration goes to Step 4 where "finalizing operations" will be performed on the set L.

*Step 3: Current iteration processes the given "current potential column selection". There are two cases.

.. Case 3.1: The given "current potential column selection" is not the complement column for a previous MajorP "column selection" – In this case, the "current potential column selection" is labelled as "column selection" for this instance of MajorP. Current iteration thereafter performs desired MajorP pivoting in the column just labelled as "column selection".

.. Case 3.2: The given "current potential column selection" is the complement column for a previous MajorP "column selection" – In this case, the iteration goes back to Step 2, to process the next not-yet-processed item in list L.

*Step 4: For each item $\hat{j}$ in L, current iteration performs MajorP pivoting in column $\hat{j}$ separately, and checks whether a solution of (Eq) has been obtained. If doing that yields a solution of (Eq), then current iteration terminates the algorithm with a solution of (Eq). Otherwise, if doing that exhausts L without yielding a solution of (Eq), then current iteration declares that (Eq) has no solutions, and terminates the algorithm at this juncture. (Claim 7.1 in Section 7 validates that conclusion)

## 4.4 A flow chart and an illustration of the iterations

We display here a summarising flow chart and a numerical illustration of the algorithm's iteration described above.

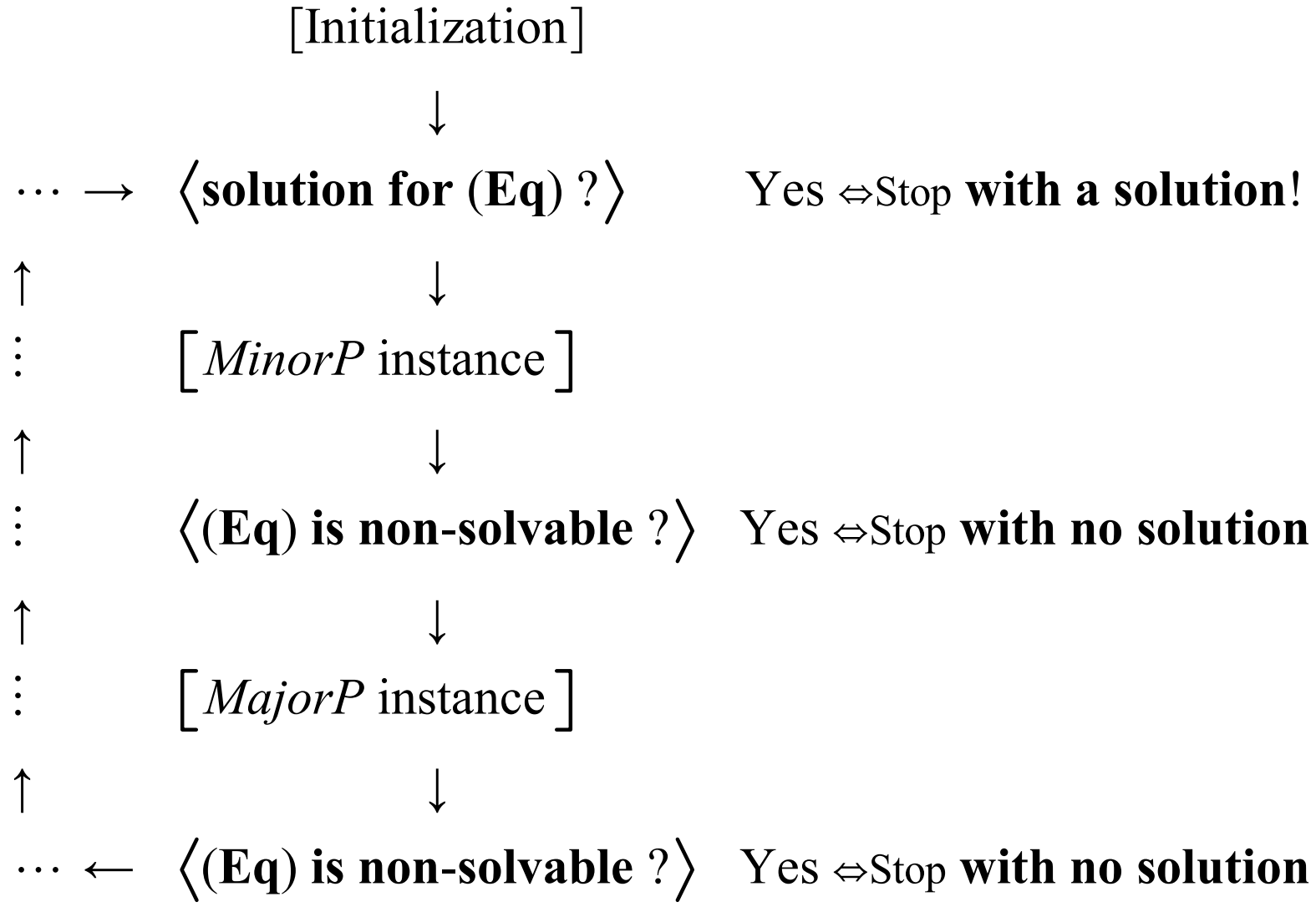

.

For a numerical illustration of the computations involved in the iterations, we will utilize the illustration LP example introduced in Section 2.

*Initialization*

$$[M\ q] = \begin{bmatrix} -10 & 5 & 0 & 2 & 1 & 0 & 0 & 0 & 10 \\ -10 & 5 & -2 & 1 & 0 & 1 & 0 & 0 & -5 \\ -11 & 6 & -1 & 1 & 0 & 0 & 1 & 0 & 1 \\ -11 & 5 & -1 & 1 & 0 & 0 & 0 & 1 & -1 \\ -10 & 5 & -1 & 1 & 0 & 0 & 0 & 0 & 0 \end{bmatrix}$$

Here, at the end of initialization, $\Pi$ is set to the empty set $\varnothing$.

.

*Iteration 1 MinorP*

L={2,4}. $\Pi = \varnothing$. Accordingly, the first instance of MinorP has its pivot at position (2,2), of current $[M\ q]$ instance to obtain the next $[M\ q]$ instance which we denote as $Z1$.

$$
Z1 = \begin{bmatrix}
0 & 0 & 2 & 1 & 1 & -1 & 0 & 0 & 15 \\
-2 & 1 & -0.4 & 0.2 & 0 & 0.2 & 0 & 0 & -1 \\
1 & 0 & 1.4 & -0.2 & 0 & -1.2 & 1 & 0 & 7 \\
-1 & 0 & 1 & 0 & 0 & -1 & 0 & 1 & 4 \\
0 & 0 & 1 & 0 & 0 & -1 & 0 & 0 & 5
\end{bmatrix}
$$

*Iteration 1 MajorP*

L={3}. On acccount of L being a singleton, $\Pi = \{3, 7\}$. Accordingly, the first instance of MajorP has its pivot at position (3,3) of current $[M\ q]$ instance, $Z1$, to obtain the next $[M\ q]$ instance which we denote as $P1$.

$$
P1 = \begin{bmatrix}
-1.4286 & 0 & 0 & 1.2857 & 1 & 0.7143 & 0 - 1.4286 & 0 & 5 \\
-1.7143 & 1 & 0 & 0.1429 & 0 & -0.1429 & 0.2857 & 0 & 1 \\
0.7143 & 0 & 1 & -0.1429 & 0 & -0.8571 & 0.7143 & 0 & 5 \\
-1.7143 & 0 & 0 & 0.1429 & 0 & -0.1429 & -0.7143 & 1 & -1 \\
-0.7143 & 0 & 0 & 0.1429 & 0 & -0.1429 & -0.7143 & 0 & 0
\end{bmatrix}
$$

*Iteration 2 MinorP*

L={4}. On acccount of L being a singleton, $\Pi = \{3, 7, 4, 8\}$. Accordingly, the second instance of MinorP has its pivot at position (4,4) of current $[M\ q]$ instance, $P1$, to obtain the next $[M\ q]$ instance which we denote as $Z2$.

$$
Z2 = \begin{bmatrix}
14 & 0 & 0 & 0 & 1 & 2 & 5 & -9 & 14 \\
0 & 1 & 0 & 0 & 0 & 0 & 1 & -1 & 2 \\
-1 & 0 & 1 & 0 & 0 & 0 & -1 & 1 & 4 \\
-12 & 0 & 0 & 1 & 0 & -1 & -5 & 7 & -7 \\
1 & 0 & 0 & 0 & 0 & 0 & 0 & -1 & 1
\end{bmatrix}
$$

*Iteration 2 MajorP*

L={1}. On acccount of L being a singleton, $\Pi = \{3, 7, 4, 8, 1, 5\}$. Accordingly, the second instance of MajorP has its pivot at position (1,1) of current $[M\ q]$, $Z2$, to obtain the next $[M\ q]$ instance which we denote as $P2$.

$$
P2 = \begin{bmatrix}
1 & 0 & 0 & 0 & 0.0714 & 0.1429 & 0.3571 & -0.6429 & 1 \\
0 & 1 & 0 & 0 & 0 & 0 & 1 & -1 & 2 \\
0 & 0 & 1 & 0 & 0.0714 & -0.8571 & 0.3571 & 0.3571 & 5 \\
0 & 0 & 0 & 1 & 0.8571 & 0.7143 & -0.7143 & -0.7143 & 5 \\
0 & 0 & 0 & 0 & -0.0714 & -0.1429 & -0.3571 & -0.3571 & 0
\end{bmatrix}
$$

The algorithm stops here with a solution of (Eq), because $q \geq 0$ & $q_5 = 0$. Thus, the illustrative example problem introduced in Section 2 has been solved in two iterations of the algorithm; it does not require Step 4 of the procedures stated in Section 4.3.

.

## 4.5 Further clarification on the iterations

(i) In Step 1 of the 4-step procedure that describes MinorP pivoting, the ordering of index set L does not have to be "descending order" of $m_{k+n+1,j} > 0$. For example, it could be in "ascending order" of $m_{k+n+1,j} > 0$. Here is one important source of variation in the algorithm, and that presents some opportunity for utilizing special structures in special classes of LP problem classes. However, it is advisable to use the same ordering of L throughout all iterations of each run of the algorithm. In comparison to that, in Step 1 of the 4-step procedure that describes MajorP pivoting, the ordering of index set L has to be a descending order, on account of Lemma 6.1 which will be stated presently.

(ii) As a closing part of "column selection" in Case 3.1 of "MajorP pivoting instance" above, one can assume, without loss of generality, that the "column selection" (that is, the selected column) has maximal $(k+n+1)$-th component, among all the columns indicated by L. That is without loss of generality, because the selected column may be multiplied by a positive number without altering the result of performing pertinent complementary GJ pivoting on that column.

(iii) One may describe our algorithm as consisting of an "initialization" procedure, followed by "MinorP → MajorP → MinorP cycles" that terminates with *either* a MajorP pivoting instance, together with a solution of (Eq), *or else* a MajorP or a MinorP pivoting instance indicating that (Eq) has no solutions. The MinorP →MajorP → MinorP cycles are a consequence of applying pairs of complementary GJ pivotings to skew-symmetric matrices, as explained in the article [19].

# 5. A report on some illustrative LP problems

We present in this Section a brief report on how our algorithm performed on some illustrative LP problems that include an instance of Klee-Minty LP problem and a Beale LP problem.

In our description of the illustrative example problems, we utilize what we call "Z-P Records Table". A Z-P Records Table displays the indices of columns of $[M\ q]$ that are selected for MinorP pivoting (that is, when $q_{k+n+1} = 0$, Z for zero) and MajorP pivoting (that is, when $q_{k+n+1} > 0$, P for positive). For the numerical example given under "An illustration of the iterations", the Z-P Records Table is:

| itn | Z | P |
|---|---|---|
| 1 | 4 | 1* |
| 2 | 2* | 3* |

We display results for two types of ordering of (the index set) $L$ for MinorP pivoting – Type 1 refers to the ascending ordering, and Type 2 refers to the descending ordering. In each example's Z-P Records Tables, column indices that belong to Π are "asterisked" (as in the iteration illustration example just stated above).

Each example LP problem is assumed to be of the Neumann symmetric form

$$\begin{aligned} \text{maximize} \quad & f^T x \\ \text{subject to:} \quad & Ax \le b \\ & x \ge 0 \end{aligned}$$

with corresponding data given in the form

$$\left( \begin{array}{c|c} f^T & \\ \hline A & b \end{array} \right)$$

.

`Example 1: An instructive LP problem`

This LP problem and its dual LP have almost-perfectly non-degenerate solutions.

$$\left(\begin{array}{cccc|c} 2 & 7 & 6 & 4 & \\ \hline 1 & 1 & 0.83 & 0.5 & 65 \\ 1.2 & 1 & 1 & 1.2 & 96 \\ 0.5 & 0.7 & 1.2 & 0.4 & 80 \end{array}\right); \quad \begin{array}{c} \text{with primal LP solution} \\ x = (0, 5.1601, 53.2015, 31.3653)^T \\ \text{and dual solution} \\ y = (6.2147, 0.7062, 0.1130)^T \end{array}$$

.

| Type 1 Z-P Records | | |
|---|---|---|
| itn | Z | P |
| 1 | 4 | 2 |
| 2 | 7 | 11 * |
| 3 | 6 | 3 |
| 4 | 5 * | 1 |

. .

| Type 2 Z-P Records | | |
|---|---|---|
| itn | Z | P |
| 1 | 5 | 1 |
| 2 | 6 | 3 |
| 3 | 7 * | 2 |

.

`Example 2: An instance of Klee-Minty LP problem (with n=3)`

This classical LP problem and its dual have a highly degenerate solutions; it is utilized in LP literature to demonstrate how some versions of Dantzig's simplex algorithm could be led astray.

$$\left(\begin{array}{ccc|c} 100 & 10 & 1 & \\ \hline 1 & 0 & 0 & 1 \\ 20 & 1 & 0 & 100 \\ 200 & 20 & 1 & 10000 \end{array}\right); \quad \begin{array}{c} \text{with primal LP solution} \\ x = (0, 0, 10000)^T \\ \text{and dual solution} \\ y = (0, 0, 1)^T \end{array}$$

.

| Type 1 Z-P Records | | |
|---|---|---|
| itn | Z | P |
| 1 | 6 | 3 |

. .

| Type 2 Z-P Records | | |
|---|---|---|
| itn | Z | P |
| 1 | 4 | 3 |
| 2 | 2 | 5 |
| 3 | 6 | 11 * |
| 4 | 8 | 10 |

.

The $n$-variable instance of Klee-Minty LP problem is solved by a version of our algorithm in exactly one iteration, with chosen pivot columns of $[M\ q]$ then being the $(2n)$-th column (by the first MinorP instance), and the $n$-th column (by the first MajorP instance).

.

`Example 3: An instance of Beale LP problem`

This is another classical LP problem that has a very degenerate solution.

$$\left(\begin{array}{cccc|c} 0.75 & -150 & 0.02 & -6 & \\ \hline 0.25 & -60 & -0.04 & 9 & 0 \\ 0.50 & -90 & -0.02 & 3 & 0 \\ 0.00 & 0 & 1.00 & 0 & 1 \end{array}\right);$$

with primal LP solution $x = (0.04, 0, 1, 0)^T$ and dual solution $y = (0, 1.5, 0.05)^T$

.

| Type 1 Z-P Records | | |
|---|---|---|
| itn | Z | P |
| 1 | 6 | 3 |
| 2 | 4 * | 2 |

. .

| Type 2 Z-P Records | | |
|---|---|---|
| itn | Z | P |
| 1 | 4 | 2 |
| 2 | 6 * | 3 |

.

```
Example 4: Another instructive example
```

$$\left(\begin{array}{cccc|c} 2 & 1 & -1 & 1 & \\ \hline 1 & 1 & 1 & 1 & 12 \\ -1 & 0 & 1 & -1 & -8 \\ 0 & 2 & 0 & -1 & 6 \end{array}\right);$$

with primal LP solution $x = (12, 0, 0, 0)^T$ and dual solution $y = (2, 0, 0)^T$

.

| Type 1 Z-P Records | | |
|---|---|---|
| itn | Z | P |
| 1 | 5 | 3 |
| 2 | 7 | 1* |
| 3 | 4 | 14 |
| 4 | 10* | 12 |

. .

| Type 2 Z-P Records | | |
|---|---|---|
| itn | Z | P |
| 1 | 2 | 4 |
| 2 | 9 | 1 |

.

```
Example 5: A problem from p.57 of Dantzig's book [5]
```

.

$$\left(\begin{array}{ccccc|c} -2 & 1 & -3 & -7 & 5 & \\ \hline 1 & 2 & 1 & 1 & 6 & 10 \\ -2 & -3 & -4 & -1 & -2 & -4 \\ 3 & 2 & 0 & 3 & 1 & 8 \end{array}\right);$$

with primal LP solution $x = (0, 0.2857, 0, 0, 1.5714)^T$ and dual solution $y = (0.9286, 0.2857, 0)^T$

.

| Type 1 Z-P Records | | |
|---|---|---|
| itn | Z | P |
| 1 | 5 | 1 |
| 2 | 3 | 8 |
| 3 | 11 * | 2 |

. .

| Type 2 Z-P Records | | |
|---|---|---|
| itn | Z | P |
| 1 | 8 | 1 |
| 2 | 2 * | 6 |
| 3 | 5 * | 3 |
| 4 | 14 * | 7 |
| 5 | 4 | 15 * |
| 6 | 11 * | 12 |

.

Example 6: Another instructive example

.

$$\left(\begin{array}{cccc|c} 3 & 4 & 1 & 7 & \\ \hline 8 & 3 & 4 & 1 & 7 \\ 2 & 6 & 1 & 5 & 3 \\ 1 & 4 & 5 & 2 & 8 \end{array}\right);$$

with primal LP solution

$$x = (0.8421, 0, 0, 0, 2632)^T$$

and dual solution

$$y = (0.0263, 1.3947, 0)^T$$

.

| Type 1 Z-P Records | | |
|---|---|---|
| itn | Z | P |
| 1 | 6 | 3 |
| 2 | 4 | 1 |
| 3 | 5 | 2 |
| 4 | 7 | 12 |
| 5 | 10 * | 13 |

. .

| Type 2 Z-P Records | | |
|---|---|---|
| itn | Z | P |
| 1 | 7 | 2 |
| 2 | 4 * | 1 |

.

Example 7: An unbouded LP problem

$$\left(\begin{array}{ccc|c} 1 & 2 & 1.5 & \\ \hline 1 & 0 & 0 & 5 \\ 0 & 1 & 0 & 6 \\ -1 & -1 & -1 & -10 \end{array}\right);$$

with primal LP solution

n. a.

and dual solution

n. a.

.

| Type 1 Z-P Records | | |
|---|---|---|
| itn | Z | P |
| 1 | 4 | 1 * |
| 2 | 6 | – |

. .

| Type 2 Z-P Records | | |
|---|---|---|
| itn | Z | P |
| 1 | 3 | 4 |
| 2 | 1 | 5 |
| 3 | 9 | 2 * |
| 4 | 6 | – |

.

```
Example 8: Another instructive example
```

This LP problem illustrates the iteration count mentioned in our proof of Lemma 7.1

$$\left(\begin{array}{cccc|c} 2 & -1 & 4 & 3 & \\ \hline 1 & 2 & 1 & 2 & 20 \\ 3 & 1 & -1 & 2 & 18 \\ 1 & 1 & -1 & -3 & 21 \end{array}\right);$$

with primal LP solution

$x = (0, 0, 20, 0)^T$

and dual solution

$y = (4, 0, 0)^T$

.

| Type 1 Z-P Records | | |
|---|---|---|
| itn | Z | P |
| 1 | 4 | 2 |
| 2 | 7 | 1 |
| 3 | 11 | 6 |
| 4 | 4 | 3 |
| 5 | 14 | 5 |
| 6 | 11 | 10 * |
| 7 | 7 | 12 * |
| 8 | 4 | 14 * |
| 9 | 9 * | 11 |

. .

| Type 2 Z-P Records | | |
|---|---|---|
| itn | Z | P |
| 1 | 6 | 1 |

# 6. Lemma 6.1 & Corollary 1

As a summary of the foregoing description, the proposed algorithm, after initialization, follows MinorP → MajorP → MinorP cycles until it either leads to a solution of (Eq), if such a solution exists, or, else, it indicates that (Eq) has no solutions.

The remainder of this article presents a validation of the proposed algorithm through its Lemma 6.1, Corollary 1 and Lemma 7.1. Lemma 6.1 provides a "necessary condition" for a MajorP column selection to be included in the basis matrix for a solution of (Eq), whereas a related corollary,

Corollary 1, provides a "sufficient condition". Lemma 7.1 puts it all together in a computational complexity statement.

Lemma 6.1 and Corollary 1 are stated and proved in this Section. This Section also includes a short note on how the algorithm utilizes Lemma 6.1 and Corollary 1 logically. In Lemma 6.1 and Corollary 1, and the remainder of this validation article up to Lemma 7.1, we will assume that each instance of the set of indices L, introduced under "General iteration", is not a singleton. Also, in the remainder of this article, in addition to assuming that $[M\ q]$ is of the form

| $m_{1,1}$ | $m_{1,2}$ | $\cdots$ | $m_{1,2(k+n)}$ | $q_1$ |
|---|---|---|---|---|
| $\vdots$ | $\vdots$ | $\ddots$ | $\vdots$ | $\vdots$ |
| $m_{k+n,1}$ | $m_{k+n,2}$ | $\cdots$ | $m_{k+n,2(k+n)}$ | $q_{k+n}$ |
| $m_{k+n+1,1}$ | $m_{k+n+1,2}$ | $\cdots$ | $m_{k+n+1,2(k+n)}$ | $q_{k+n+1}$ |

we will assume, without any loss of generality, that each instance of $[M\ q]$ has its columns enumerated or re-enumerated such that, for $j = 1,\ldots,k+n,$ column $j$ is the complement for column $k+n+j,$ and columns $k+n+1,\ldots,2(k+n)$ form the (k+n+1)-by-(k+n) submatrix

$$U = \begin{pmatrix} 1 & 0 & 0 & & \cdots & 0 \\ 0 & 1 & 0 & & \cdots & 0 \\ 0 & 0 & 1 & & & \vdots \\ \vdots & \vdots & 0 & & \ddots & 0 \\ 0 & 0 & \cdots & 0 & & 1 \\ 0 & 0 & \cdots & 0 & & 0 \end{pmatrix}$$

# 6.1 Lemma 6.1 and its proof

Lemma 6.1 is about MajorP column selections. Towards stating Lemma 6.1, we begin by informally stating some 'background clarification notes' that outline a context that may aid a clear understanding of Lemma 6.1 and our proof of it.

## 6.1.1 Background clarification notes.

(i) – When, in Section 4.3, a MajorP instance selects a maximal column of $[M\ q],$ say column $M^{(2)}$ for example, that does not mean that column $M^{(2)}$ will certainly be included in the basis matrix for some solution of (Eq). Instead, as we will show shortly through Lemma 6.1, such a selection of $M^{(2)}$ only means that $M^{(2)}$ satisfies a *necessary condition* for $M^{(2)}$ to be included in the basis matrix for a solution of (Eq). Corollary 1 is about a related *sufficient condition.*

(ii) – Each column that is selected by a MajorP instance has a positive (k+n+1)-th component, and the selected column is only 'competing' against its own complement column, which is always a unit-vector in current $[M\ q]$ instance.

(iii) – For a MajorP pivoting, $q_{k+n+1} > 0.$ If (Eq) has any solutions, then $q_{k+n+1} > 0$ 'forces' at least one of the columns of $M$ having positive (k+n+1)-th component to be utilized in some solution of (Eq). In Lemma 6.1 statement below, such a "forced column" works hand-in-hand with a "maximal column" (that is a column of $M$ that has a maximal (k+n+1)-th component) to take one step forward towards building a basis matrix for a solution of (Eq).

(iv) – The last row of each $[M\ q]$, which is usually regarded as redundant in solving LP problems as Linear Complementarity Problems, turns out to be an indispensable extra dimension in our proof of Lemma 6.1, especially in obtaining 'arbitrary closeness to 0' in the proof.

.

LEMMA 6.1: *Suppose that the first column* $M^{(1)}$ *and the second column* $M^{(2)}$ *of current* $[M\ q]$ *are to be compared for a MajorP pivoting. Suppose that* $m_{k+n+1,2} > m_{k+n+1,1} > 0$, *in row* k + n + 1 *of* $[M\ q]$.

*If there is a solution of (Eq), say* $z^*$, *utilizing column* $M^{(1)}$ *and/or column* $M^{(2)}$ *such that* $z_1^* + z_2^* > 0$, *then there is a solution* $\bar{z}$ *of the sub-system* $Mz = q$ *of (Eq), such that* $\bar{z}_1 = 0$ *and* $\bar{z}_2 > 0$ (*that is, utilizing* $M^{(2)}$ *but not* $M^{(1)}$). *Moreover, for* $j = 1,\ldots,2(k+n)$, $\bar{z}_j > 0$ *where* $z_j^* > 0$, *and* $\bar{z}_j$ *is 0 or is 'arbitrarily close to* $0'$ *where* $z_j^* = 0$, *with that 'arbitrarily close to* 0' *effected by adding a suitably large multiple of row* $k+n+1$ *of* $[M\ q]$ *to row i of* $[M\ q]$, *for* $i = 1,\ldots,k+n$ *as needed.*

.

Before stating a proof of Lemma 6.1, we observe here, for an emphasis, that $\bar{z}$ in the conclusion of Lemma 6.1 would have constituted a solution of (Eq), instead of just being a solution of the sub-system $Mz = q$ of (Eq), if the hypothesis of Lemma 6.1 had included (i) the statement $z_{k+n+2}^* = 0$ and (ii) "arbitrarily close to 0" is equated to 0. That is true because $\bar{z}_{k+n+2}$ would then be "arbitrarily close to 0" as well. This observation is intended to indicate a plausibility about the conclusion of Lemma 6.1, and it may also be regarded as a pointer to Corollary 1 below.

.

PROOF OF LEMMA 6.1

.

If $z_1^* = 0$, then the conclusion of Lemma 6.1 is obtained trivially. So hencforth assume $z_1^* > 0$. The strategy of this proof is to *express* $Mz^* = q$ suitably (in the vector equation (2) below) and thereafter *propose* $M\bar{z} = q$ (in vector equation (3) below), in a form that resembles (2), in order to enable a comparison of the coefficients of pertinent columns $M^{(j)}$'s in each one of the k + n + 1 component equations that (2) and (3) become.

We will assume that the components of the two columns $M^{(1)}$ & $M^{(2)}$ are positive, since we could add a suitable multiple of row $k+n+1$ of $[M\ q]$ to the other rows to make it true.

From $M\,z^* = q$, we have

$$z_1^* M^{(1)} + z_2^* M^{(2)} + \sum_{j=3}^{j=2(k+n)} z_j^* . M^{(j)} = q \cdots\cdots(1)$$

Assuming that columns $k+n+1,\ldots,2(k+n)$ of $M$ comprise the matrix $U \equiv [u^{(1)}|\cdots|u^{(k+n)}]$ defined earlier, re-write (1) as

.

$$z_1^* M^{(1)} + \sum_{j=3}^{j=2(k+n)} z_j^* . M^{(j)} = q - z_2^* M^{(2)}$$

that is,

$$z_1^* M^{(1)} + \sum_{j=k+n+1}^{j=2(k+n)} z_j^* . u^{(j-k-n)} = q - z_2^* M^{(2)} - \sum_{j=3}^{j=k+n} z_j^* . M^{(j)}$$

$$= r \qquad \cdots\cdots(2)$$

where, $r \equiv q - z_2^* M^{(2)} - \sum_{j=3}^{j=k+n} z_j^* . M^{(j)}$, and $r_{k+n+1} > 0$ on acccount of $m_{k+n+1,1} > 0, z_1^* > 0$ and $z^*$ being a solution of (Eq).

That is, in terms of its component equations,

$$z_1^* m_{i,1} + z_{k+n+i}^* = r_i, \text{ for } i = 1,\ldots,k+n; \qquad \cdots\cdots(2)$$
$$\& \; z_1^* m_{k+n+1,1} = r_{k+n+1}$$

Towards the conclusion of Lemma 6.1, we begin by setting

$$\boxed{\begin{array}{l} \bar{z}_1 = 0, \text{ and} \\ \bar{z}_j = z_j^* \text{ for } j = 3,\ldots,k+n, \end{array}}$$

For the purpose of determining the $\bar{z}_j$'s for $j = k+n+1,\ldots,2(k+n)$ in the conclusion of Lemma 6.1, there is a good deal of analysis work to be done. We will, treating $M^{(1)}$ & $M^{(2)}$ with the same basic matrix operations in $[M\,q]$ as context, begin by staying close to the structure of (2), and propose the vector equation

$$\bar{z}_2 M^{(2)} + \sum_{j=k+n+1}^{j=2(k+n)} \bar{z}_j . u^{(j-k-n)} = r \cdots\cdots(3)$$

Note that the right-hand-side $r$ is the same for both (2) and (3). In terms of its component equations, in a manner similar to (2), the vector equation (3) becomes

$$\bar{z}_2 m_{i,2} + \bar{z}_{k+n+i} = r_i, \text{ for } i = 1,\ldots,k+n; \qquad \cdots\cdots(3)$$
$$\& \; \bar{z}_2 m_{k+n+1,2} = r_{k+n+1}$$

Note how we use index $i$ for rows of $M$ and use index $j$ for the columns of $M$. The remainder of this proof consists of determining $\bar{z}_2$ & $\bar{z}_j, j = k+n+1,\ldots,2(k+n)$, such that equations (3) hold, along with

$\bar{z}_2 > 0$; $\bar{z}_j > 0$ where $z_j^* > 0$ and $\bar{z}_j$ is 0 or is "arbitrarily close to 0" where $z_j^* = 0$,

for $j = k+n+1,\ldots,2(k+n)$, with "arbitrarily close to 0" effected through (the elementary row operation of) adding a suitably large multiple of the last row of $[M\,q]$ to some row of $[M\,q]$.

.

Towards that task of determining the $\bar{z}_j$'s, we will utilize a comparison of the component/row equations comprising (2) and the component/row equations comprising (3), in five distinct cases, namely, Case 0, Case 1, Case 2, Case 3, Case 4 below.

.

`Case 0: Component/row` $k+n+1$ `equation` (i.e. last row of $[M\,q]$)

.

From (2), $z_1^* m_{k+n+1,1} = r_{k+n+1}$. From (3), $\bar{z}_2 m_{k+n+1,2} = r_{k+n+1}$. Therefore, to determine $\bar{z}_2$ we set $z_1^* m_{k+n+1,1} = \bar{z}_2 m_{k+n+1,2}$; that is, with $f \equiv m_{k+n+1,1}/m_{k+n+1,2}$,

$$\boxed{\bar{z}_2 = z_1^* f, \text{ for Case 0}}$$

Thus, $\bar{z}_2 > 0$ as $z_1^* > 0$.

.

`Case 1: Component/row` $i$ `equation, with` $z_{k+n+i}^* > 0$ & $m_{i,1} \geq m_{i,2} > 0$ (in an i-th row of $[M\,q]$)

.

Because $z_1^* > 0$, $z_1^* m_{i,1} \geq z_1^* m_{i,2} > 0$. Because $f < 1$, there is a positive number, say $\mu_i$, such that $z_1^* m_{i,1} = z_1^* m_{i,2}.f + \mu_i$. That is,

$$\begin{aligned} z_1^* m_{i,1} &= z_1^* m_{i,2}.f + \mu_i. \\ &= \bar{z}_2 m_{i,2} + \mu_i \end{aligned}$$

That is, from (2),

$$r_i = z_1^* m_{i,1} + z_{k+n+i}^* = \bar{z}_2 m_{i,2} + \mu_i + z_{k+n+i}^*$$

In (3), we then set

$$\boxed{\bar{z}_{k+n+i} = z_{k+n+i}^* + \mu_i > 0 \text{ for Case 1}}$$

For Cases 2 & 3 below, we define

$$\theta_i = (m_{i,1} - m_{i,2})/(m_{k+n+1,2} - m_{k+n+1,1})$$

Considering the elementary row operation of adding $\theta_i$ multiple of the last row (row $k+n+1$) of $[M\, q]$ to its row $i$, we will sometimes refer to $\theta_i$ as "component $i$ equalizer for columns $M^{(1)}$ & $M^{(2)}$", since

$$m_{i,1} + \theta_i . m_{k+n+1,1} = m_{i,2} + \theta_i . m_{k+n+1,2}$$

.

An outline of what we will do in Cases 2, 3 & 4 is as follows. We will consider whether $z_1^* . (m_{i,1} + \theta_i . m_{k+n+1,1}) + z_{k+n+i}^*$ is negative or not for each $i$ in Cases 2 &3. If $z_1^* . (m_{i,1} + \theta_i . m_{k+n+1,1}) + z_{k+n+i}^* \geq 0$ (Case 2), then a straightforward comparison of

$$z_1^* . (m_{i,1} + \theta_i . m_{k+n+1,1}) + z_{k+n+i}^* \text{ \& } [z_1^* . (\mathrm{m}_{i,2} + \theta_i . \mathrm{m}_{k+n+1,2}) + z_{k+n+i}^*] f$$

.

will yield the desired conclusion. If $z_1^* . (m_{i,1} + \theta_i . m_{k+n+1,1}) + z_{k+n+i}^* < 0$ (Case 3), then increasing the magnitude of $\theta_i$ enough will lead to the desired conclusion after the intermediate value theorem is utilized. In the case of $z_{k+n+i}^* = 0$ (Case 4), a normalization of row $i$ will lead to the desired conclusion.

.

Case 2: Component/row $i$ equation, with $z_{k+n+i}^* > 0$; $m_{i,1} < m_{i,2}$ & $z_1^* . (m_{i,1} + \theta_i . m_{k+n+1,1}) + z_{k+n+i}^* \geq 0$ (in an i-th row of $[M\, q]$)

.

From $m_{i,1} + \theta_i . m_{k+n+1,1} = m_{i,2} + \theta_i . m_{k+n+1,2}$ and $z_1^* . (m_{i,1} + \theta_i . m_{k+n+1,1}) + z_{k+n+i}^* \geq 0$, we have

$$z_1^* . (\mathrm{m}_{i,1} + \theta_i . \mathrm{m}_{k+n+1,1}) + z_{k+n+i}^* = z_1^* . (\mathrm{m}_{i,2} + \theta_i . \mathrm{m}_{k+n+1,2}) + z_{k+n+i}^* \geq 0.$$

Because $f < 1$, one can then find a non-negative number, say $\delta_i$, such that

$$z_1^* . (\mathrm{m}_{i,1} + \theta_i . \mathrm{m}_{k+n+1,1}) + z_{k+n+i}^* = [z_1^* . (\mathrm{m}_{i,2} + \theta_i . \mathrm{m}_{k+n+1,2}) + z_{k+n+i}^*] f + \delta_i.$$

From $\theta_i . \mathrm{m}_{k+n+1,1} = \theta_i . \mathrm{m}_{k+n+1,2} . f$, on account of the definition of $f$, that becomes

$$z_1^* . \mathrm{m}_{i,1} + z_{k+n+i}^* = z_1^* . \mathrm{m}_{i,2} . f + z_{k+n+i}^* . f + \delta_i,$$

That is, using $\bar{z}_2 = z_1^* f$,

$$z_1^* . \mathrm{m}_{i,1} + z_{k+n+i}^* = \bar{z}_2 . \mathrm{m}_{i,2} + z_{k+n+i}^* . f + \delta_i,$$

Therefore, from comparing (2) and (3), we set

$$\boxed{\bar{z}_{k+n+i} = z_{k+n+i}^* f + \delta_i > 0 \text{ for Case 2}}$$

Case 3: Component/row $i$ equation, with $z_{k+n+i}^* > 0$; $m_{i,1} < m_{i,2}$; & $z_1^* . (m_{i,1} + \theta_i . m_{k+n+1,1}) + z_{k+n+i}^* < 0$ (in an i-th row of $[M\, q]$)

.

Define $\widetilde{\theta_i}$ by $\widetilde{\theta_i} \equiv \theta_i - \varepsilon_i$, where $\varepsilon_i$ is a positive number that one may choose close to 0. Then $\widetilde{\theta_i} < \theta_i = (\mathrm{m}_{i,1} - \mathrm{m}_{i,2})/(\mathrm{m}_{k+n+1,2} - \mathrm{m}_{k+n+1,1})$, and $z_1^* . (m_{i,1} + \widetilde{\theta}_i . m_{k+n+1,1}) + z_{k+n+i}^* < 0$ on account of $z_1^* . (m_{i,1} + \theta_i . m_{k+n+1,1}) + z_{k+n+i}^* < 0$ and $m_{k+n+1,1} > 0$.

Using $(\mathrm{m}_{i,2} + \widetilde{\theta}_i . \mathrm{m}_{k+n+1,2}) < (\mathrm{m}_{i,1} + \widetilde{\theta}_i . \mathrm{m}_{k+n+1,1})$ (from $\widetilde{\theta_i}$'s definition), and $z_1^* . (m_{i,1} + \widetilde{\theta}_i . m_{k+n+1,1}) + z_{k+n+i}^* < 0$ (from this case's hypothesis), we have

$$0 > z_1^* . (\mathrm{m}_{i,1} + \widetilde{\theta}_i . \mathrm{m}_{k+n+1,1}) + z_{k+n+i}^* > z_1^* . (\mathrm{m}_{i,2} + \widetilde{\theta}_i . \mathrm{m}_{k+n+1,2}) + z_{k+n+i}^*$$

Define a function $H : [0,1] \to R$ by

$$H(v) \equiv z_1^* . (\mathrm{m}_{i,1} + \widetilde{\theta}_i . \mathrm{m}_{k+n+1,1}) + z_{k+n+i}^* - [z_1^* . (\mathrm{m}_{i,2} + \widetilde{\theta}_i . \mathrm{m}_{k+n+1,2}) + z_{k+n+i}^*] . v \cdots\cdots (4)$$

Clearly $H(0) < 0 < H(1)$. Utilizing an intermediate value theorem, one can find a positive number $v^*$ in the interval $(0,1)$ such that $H(v^*) = 0$; that is, $v^*$ a zero of $H$, and

$$z_1^*.(m_{i,1} + \widetilde{\theta}_i.m_{k+n+1,1}) + z_{k+n+i}^* = [z_1^*.(m_{i,2} + \widetilde{\theta}_i.m_{k+n+1,2}) + z_{k+n+i}^*].v^*$$

The number $v^*$ depends on $i$, and, to reflect that dependence, we will write $v^*$ as $v_i^*$. Towards desired conclusion, there are two sub-cases to consider at this juncture – the sub-case $v_i^* \leq f$ and the sub-case $v_i^* > f$.

.

*Sub-case 3.1:* $v_i^* \leq f$

For $f$ in the interval $[v_i^*, 1)$, one can see that

$$[z_1^*.(m_{i,2} + \widetilde{\theta}_i.m_{k+n+1,2}) + z_{k+n+i}^*].v_i^* \geq [z_1^*.(m_{i,2} + \widetilde{\theta}_i.m_{k+n+1,2}) + z_{k+n+i}^*].f,$$

because $z_1^*.(m_{i,2} + \widetilde{\theta}_i.m_{k+n+1,2}) + z_{k+n+i}^* < 0$. Therefore, for $f$ in the interval $[v_i^*, 1)$, one can see that

$$\begin{aligned} z_1^*.(m_{i,1} + \widetilde{\theta}_i.m_{k+n+1,1}) + z_{k+n+i}^* &= [z_1^*.(m_{i,2} + \widetilde{\theta}_i.m_{k+n+1,2}) + z_{k+n+i}^*].v_i^* \\ &\qquad \text{(because } H(v_i^*) = 0\text{)} \\ &\geq [z_1^*.(m_{i,2} + \widetilde{\theta}_i.m_{k+n+1,2}) + z_{k+n+i}^*].f; \end{aligned}$$

that is,

$$z_1^*.m_{i,1} + z_{k+n+i}^* \geq z_1^*.m_{i,2}.f + z_{k+n+i}^*.f,$$

on account of the definition $f \equiv m_{k+n+1,1}/m_{k+n+1,2}$. One can then find a number $\lambda_i \geq 0$ such that $z_1^*.m_{i,1} + z_{k+n+i}^* = z_1^*.m_{i,2}.f + z_{k+n+i}^* f + \lambda_i$;

that is, from using $\bar{z}_2 = z_1^* f$,

$$z_1^*.m_{i,1} + z_{k+n+i}^* = \bar{z}_2.m_{i,2}. + z_{k+n+i}^* f + \lambda_i.$$

Accordingly, in (3), we set

$$\boxed{\bar{z}_{k+n+i} = z_{k+n+i}^*.f + \lambda_i > 0 \text{ for } \textit{Sub-case 3.1}}$$

That is not contingent on the value of $\widetilde{\theta}_i$, because $\widetilde{\theta}_i$ is cancelled out by $f \equiv m_{k+n+1,1}/m_{k+n+1,2}$.

.

*Sub-case 3.2:* $v_i^* > f$

Recall from (4) that

$$H(v) \equiv z_1^*.(m_{i,1} + \widetilde{\theta}_i.m_{k+n+1,1}) + z_{k+n+i}^* - [z_1^*.(m_{i,2} + \widetilde{\theta}_i.m_{k+n+1,2}) + z_{k+n+i}^*].v$$

One can decrease the value of $v_i^*$ by increasing $|\widetilde{\theta}_i|$, because $m_{k+n+1,2} > m_{k+n+1,1}$ and $H(v)$ increases for every $v$ as $|\widetilde{\theta}_i|$ is increased. A sufficiently large increase of $|\widetilde{\theta}_i|$ would result in $v_i^* < f$. Thus, *Sub-case 3.2* reduces to an instance of *Sub-case 3.1*. Accordingly, in (3), we set

$$\boxed{\bar{z}_{k+n+i} = z_{k+n+i}^*.f + \lambda_i > 0 \text{ for } \textit{Sub-case 3.2} \text{ as well}}$$

Since there is only a finite number of $i$'s that are in Sub-case 3.2, one can find one value of $\widetilde{\theta}_i$ that will work for all such $i$'s.

.

Case 4 – All components/rows wherein $z_{k+n+i}^* = 0$

.

Here we will suitably normalize row $i$, and we begin with two needed elementary row operations on row $i$ of $[M\ q]$ (that is, equation $i$ of (2)) having $z_{k+n+i}^* = 0$. Choose a large positive number, say $\beta$, and add $\beta$ multiple of equation $k + n + 1$ of (2) to equation $i$ of (2) having $z_{k+n+i}^* = 0$; thereafter, divide both sides of resultant equation $i$ of (2) by its right-hand-side $r_i + \beta.r_{k+n+1}$.

Accordingly, with $z_{k+n+i}^* = 0$ and $i \leq k + n$, the i-th equation of (2) becomes (2i).

$$z_1^*.(\mathrm{m}_{i,1} + \beta.\mathrm{m}_{k+n+1,1})/(r_i + \beta.r_{k+n+1}) = 1 \cdots\cdots(2i)$$

To determine $\bar{z}_{k+n+i}$ that corresponds to $z_{k+n+i}^* = 0$ in equation $(2i)$, we propose $(3i)$ in the same vein that (3) was proposed:

$$\bar{z}_2(m_{i,2} + \beta.\mathrm{m}_{k+n+1,2})/(r_i + \beta.r_{k+n+1}) + \bar{z}_{k+n+i} = 1 \cdots\cdots(3i)$$

From $(2i)$ and $(3i)$

$$\begin{aligned}
\bar{z}_{k+n+i} &= [z_1^*.(\mathrm{m}_{i,1}+\beta.\mathrm{m}_{k+n+1,1}) - \bar{z}_2(\mathrm{m}_{i,2}+\beta.\mathrm{m}_{k+n+1,2})]/(\mathrm{r}_i+\beta.\mathrm{r}_{k+n+1})] \\
&= [z_1^*.(\mathrm{m}_{i,1}/\beta + \mathrm{m}_{k+n+1,1}) - \bar{z}_2(\mathrm{m}_{i,2}/\beta + \mathrm{m}_{k+n+1,2})]/(\mathrm{r}_i/\beta + \mathrm{r}_{k+n+1})] \\
&= z_1^*.[(\mathrm{m}_{i,1}/\beta + \mathrm{m}_{k+n+1,1}) - f.(\mathrm{m}_{i,2}/\beta + \mathrm{m}_{k+n+1,2})]/(\mathrm{r}_i/\beta + \mathrm{r}_{k+n+1})] \\
&= z_1^*.[(\mathrm{m}_{i,1}/\beta) - f.(\mathrm{m}_{i,2}/\beta)]/(\mathrm{r}_i/\beta + \mathrm{r}_{k+n+1})]
\end{aligned}$$

Note that each one of the four numbers

$$\mathrm{m}_{i,1}/\beta,\ f.\mathrm{m}_{i,2}/\beta,\ [(\mathrm{m}_{i,1}/\beta) - f.(\mathrm{m}_{i,2}/\beta)]\ \&\ \mathrm{r}_i/\beta$$

becomes very small as β is chosen to be very large. That is, $\bar{z}_{k+n+i}$ gets arbitrarily close to $z_1^*[0-0]/(0+\mathrm{r}_{k+n+1})$ as β is chosen arbitrarily large. Thus, one can choose β arbitrarily large to have $\bar{z}_{k+n+i}$ arbitrarily close to 0. Note also that the choice of β may be considered not contingent upon $i$, because there is only a finite number of $i$'s having $z_{k+n+i}^* = 0$.

.

END OF PROOF OF LEMMA 6.1

.

In the algorithm, "arbitrarily close to 0" (Case 4) is automatically equated to "0" through complementary GJ pivoting maintaining $z_j z_{(k+n+j)} = 0$, for $j = 1, \ldots, k+n$ throughout.

.

## 6.2 Corollary 1 and its proof

We next state and prove Corollary 1 which is essentially a corollary of Lemma 6.1. Corollary 1 will demonstrate two things:

(a) that $\bar{z}$, the solution of sub-system $Mz = q$ (in the conclusion of Lemma 6.1) becomes a solution of (Eq), if the hypothesis of Lemma 6.1 includes the statement $z_{k+n+2}^* = 0$ and the phrase "arbitrarily close to 0" (in the conclusion of Lemma 6.1) is equated to "0"; and

(b) the algorithm's MajorP pivoting on $M^{(2)}$ (in the hypothesis of Lemma 6.1) would detect a solution of (Eq), and that solution will turn out to be $\bar{z}$.

Before stating Corollary 1, we first state some "Background clarification notes" that are intended to aid a clear understanding of Corollary 1.

.

### 6.2.1 Background clarification notes for Corollary 1

(i) – Recall that we are assuming that each $[M\ q]$ instance has its columns enumerated or re-enumerated such that, for $\mathrm{j} = 1, \ldots, \mathrm{k+n}$, column $j$ is the complement for column $k+n+j$, and that columns $\mathrm{k+n+1}, \ldots, 2(\mathrm{k+n})$ form the (k+n+1)-by-(k+n) sub-matrix U. The $k+n$ unit vectors that comprise U correspond to some k+n columns of the very first $[M\ q]$ introduced at the algorithm's "Initialization". Henceforth, we will refer to that first $[M\ q]$ as "the initial $[M\ q]$".

(ii) If $q_{k+n+1} = 0$ in current $[M\ q]$ instance, then, as a definition here, the *k+n columns of the initial* $[M\ q]$ that correspond to the U submatrix of current $[M\ q]$ constitute a *basis matrix* for a solution of (Eq), a solution that may be feasible or infeasible; it is infeasible if there is a $j$ such that $q_j < 0$ in current $[M\ q]$ instance.

(iii) Complementary GJ pivoting in the algorithm ensures that the complementary slackness

requirements, $z_j z_{(k+n+j)} = 0,$ for $j = 1,\ldots,k+n$ (stated in (Eq)) are satisfied throughout. This effect of complementary GJ pivoting in the algorithm is what enables one to interpret "arbitrarily close to 0" (in the conclusion of Lemma 6.1) as "0", whenever such an interpretation is needed.

(iv) Each $\bar{z}_j$ that is "arbitrarily close to 0" may be interpreted as "0" in the conclusion of Lemma 6.1 because corresponding $\beta$ (specified in our proof of Lemma 6.1) is free to be made as large as needed, without affecting any other components of $\bar{z}$ or of the algorithm.

.

COROLLARY 1: *In Lemma 6.1, suppose the hypothesis also says that $z^*_{k+n+2} = 0$ in the given solution $z^*$ of (Eq). Under that condition:*

*(a) $\bar{z}_{k+n+2}$ too is "arbitrarily close to 0" in the conclusion of Lemma 6.1*;

*(b) if each "arbitrarily close to 0" component of it is set equal to 0, then $\bar{z}$ in the conclusion of Lemma 6.1 is a feasible solution of (Eq)*;

*(c) the algorithm's complementary MajorP pivoting in $M^{(2)}$ (of Lemma 6.1 statement) would produce a feasible solution of (Eq) that is exactly $\bar{z}$. That feasible solution is $(B^TB)^{-1}B^Tq$ suitably filled-up with zeros, where $B$ is basis matrix (corresponding to $\bar{z}$) inside the initial $[M\,q]$* .

.

PROOF OF COROLLARY 1

.

Recall that we continue to assume that columns of $[M\,q]$ are enumerated such that columns $k+n+1,\ldots,2(k+n)$ of $[M\,q]$ are unit vectors $u^{(1)},\ldots,u^{(k+n)}$ respectively in $R^{k+n+1}$, and are complements for columns $1,\ldots,k+n$ respectively.

(a) With $z^*_{k+n+2} = 0$ in the hypothesis of Lemma 6.1, $\bar{z}_{k+n+2}$ too would be "arbitrarily close to 0" in the conclusion of Lemma 6.1, as already demonstrated in the proof of Lemma 6.1.

(b) Since each $\bar{z}_i$ that is "arbitrarily close to 0" is now set to "0" without changing the other positive $\bar{z}_j$'s, the vector $\bar{z}$ satisfies not only $Mz = q$ as in the conclusion of Lemma 6.1, but also $z \geq 0$ and $z_j z_{(k+n+j)} = 0,$ for $j = 1,\ldots,k+n.$ Therefore, $\bar{z}$ is a feasible solution of (Eq).

(c) Consider MajorP pivoting in column $M^{(2)}$ of current $[M\,q]$ (as in Lemma 6.1 statement). That MajorP pivoting has the effect, in the initial $[M\,q]$, of using the column corresponding to $M^{(2)}$ to displace the column corresponding to $M^{(k+n+2)}$ from a basis matrix, thereby resulting in a new basis matrix, say $B$. Next recall how $\bar{z}$ is defined in Lemma 6.1. On account of $\bar{z}_j > 0$ for every $z^*_j > 0,$ the basis matrix corresponding to $\bar{z}$ in the initial $[M\,q]$ only exchanges the column corresponding to $M^{(1)}$ with the column corresponding to $M^{(2)}$, thereby resulting in the same basis matrix $B$ in the initial $[M\,q]$. From basic linear algebra, the system of equations $Bx = q,$ has at most one solution $x$, which is $(B^TB)^{-1}B^Tq$; clearly, it has one solution through $z.$ It follows that the 2(k+n)-vector $\bar{z}$ is the (k+n)-vector $(B^TB)^{-1}B^Tq$ suitably extended with zeros. Thus, MajorP pivoting in $M^{(2)}$ would detect $\bar{z}$ as a feasible solution of (Eq).

The following diagram depicts the reasoning just stated in (c) above.

.

| $\bar{z}$ from Lemma 6.1 and Corollary 1 | $\longrightarrow$ | ( the same basis matrix, $B$, in the initial $[Mq]$ ) | $\longleftarrow$ | MajorP pivoting in $M^{(2)}$of current $[Mq]$ |
|---|---|---|---|---|

$$\downarrow$$

| only solution |
|---|
| $x = (B^T B)^{-1} B^T q$ |
| of $Bx = q$, |
| so that $\bar{z} \equiv x$ |

.

END OF PROOF OF COROLLARY 1

.

## 6.3 On how the algorithm utilizes Lemma 6.1 and Corollary 1

Lemma 6.1 and Corollary 1 are about MajorP column selections. When a column selection is made by MajorP, it means that the column has passed a test that is 'necessary' for that column to be included in the basis matrix for a solution of (Eq). A corresponding 'sufficient' condition test is provided by Corollary 1.

Lemma 6.1's role may be regarded as that of gradually assembling columns of initial $[M\ q]$ that may form a basis matrix for a feasible solution of (Eq). Corollary 1 is about when local information (e.g. $z^*_{k+n+2} = 0$) indicates that a basis matrix for a feasible solution of (Eq) is available.

We will demonstrate in the next Section that such local information exists when the MinorP instance → MajorP instance → MinorP instance → $\cdots$ process unavoidably reverses certain previous MajorP column selections, unless (Eq) has no solutions.

Details of that last statement are provided through Claim 7.1, Claim 7.2 and Lemma 7.1 under "Computational complexity statement and its proof".

# 7. Computational complexity statement and its proof

In this Section, we will explain the computational efficiency/complexity of the algorithm. We will do this by tracking what happens to MajorP column selections in the algorithm. To do that tracking, we will consider two cases, namely, the case wherein a previous MajorP column selection gets reversed along the way, and the case wherein no previous MajorP selections get reversed along the way.

Towards that, in the "Preliminary" section below, we will show (through two Claims – Claims 7.1 & 7.2) that, if (Eq) has a solution, then a MajorP selection reversal causes the hypotheses of Lemma 6.1 and Corollary 1 to hold. After the "Preliminary", we will state and prove our computational complexity lemma, Lemma 7.1.

## 7.1 Preliminary for computational complexity statement

Each one of Claims 7.1 & 7.2 is about local information becoming available for satisfying the hypotheses of Lemma 6.1 and Corollary 1, whenever a previous MajorP column selection *has to be reversed and* (Eq) *has a solution.*

Regarding notation for Claim 7.1 and Claim 7.2, $[M\ q]$ will refer to the current $[M\ q]$ instance, even though the statements are true for all $[M\ q]$ instances as well.

### 7.1.1 When a MajorP instance reverses a previous MajorP selection

Recall that a MajorP instance can reverse a previous MajorP selection only in Step 1 and Step 4 of MajorP instance description of Section 4.3 (Executing next iteration).

The reversal in Step 1 is on account of L being a singleton; that translates into an implicit reduction, as the corresponding column of $M$ is implicitly put aside as "settled" or "eliminated".

In the reversal at Step 4 of Section 4.3, each item in list L is the index for the complement of a previous MajorP column selection. As the list L contains the indices of all other columns of $[M\ q]$ that have positive (k+n+1)-th component, one can see that $q_{k+n+1} > 0$ forces at least one of the columns indexed by L to be included in the basis matrix of a solution of (Eq), if (Eq) has any solutions. Claim 7.1 is about the result of a reversal (of a previous MajorP selection) of this type.

.

CLAIM 7.1: *Consider the algorithm processing Step 4 of "MajorP pivoting instance" in Section 4.3. If (Eq) has a solution, then there is a column of* $[M\ q]$, *denoted here as* $M^{(t)}$ *with index t contained in list L, such that the complementary GJ pivoting in* $M^{(t)}$ *terminates the algorithm with a solution of (Eq), as specified in Corollary 1.*

.

PROOF OF CLAIM 7.1

.

Recall that the list $L$, at this juncture, has at least two elements, say $v$ & $w$. Let $M^{(v)}$ be maximal in row k+n+1, and let $M^{(w)}$ be a column utilized in the basis matrix for a solution, say $z^*$, of (Eq). For notation convenience, assume that $v \leq k+n$ & $w \leq k+n$, so that $M^{(k+n+v)}$ and $M^{(k+n+w)}$ are complements of $M^{(v)}$ and $M^{(w)}$ respectively. Also, note that some previous MajorP instances did their complementary GJ pivotings in $M^{(k+n+v)}$ and $M^{(k+n+w)}$ (which are $v$-th unit vector and $w$-th unit vector respectively, in current $[M\ q]$). The following diagram is intended to aid some intuition while introducing iterations $i_1$, $i_2$ & $i_3$ as well.

| $M^{(k+n+v)}$ at iteration $i_1$ | & | $M^{(k+n+w)}$ at iteration $i_2$ | , … , | $M^{(v)}$ & $M^{(w)}$ at iteration $i_3$, with $i_1, i_2 < i_3$ |
|---|---|---|---|---|

We consider two cases, the case $v \neq w$ and the case $v = w$.

Case 1: $v \neq w$

By virtue of Lemma 6.1, as $z^*_v + z^*_w > 0$, there is a solution $\bar{z}$ of the sub-system $Mz = q$ of (Eq), such that $\bar{z}_w = 0$ and $\bar{z}_v > 0$ (that is, utilizing $M^{(v)}$ but not $M^{(w)}$). By virtue of an LP complementary slackness theorem, there are two sub-cases to consider here, namely, the subcase $z^*_v \neq 0$ and the subcase $z^*_v = 0$.

Subcase 1.1: $z^*_v \neq 0$. By an LP complementary slackness theorem, $z^*_{k+n+v} = 0$ and $z_{k+n+v} = 0$ for every solution $z$ of (Eq). By Lemma 6.1, $\bar{z}_{k+n+v}$ is arbitrarily close to 0. By Corollary 1, a complementary MajorP pivoting in $M^{(v)}$ would terminate the algorithm with a solution of (Eq). Therefore, in this subcase, Subcase 1.1, we set $t \leftarrow v$ (t introduced in Claim 7.1).

Subcase 1.2: $z^*_v = 0$. By virtue of an LP strict complementarity theorem, $z^*_v = 0$ for every solution $z^*$ of (Eq) (see Chapter 3 of [14] for a straightforward statement and proof of that LP strict complementarity theorem). And by virtue of Lemma 6.1 & Corollary 1, the earlier complementary MajorP pivoting in $M^{(k+n+v)}$ (at iteration $i_1$ wherein $M^{(k+n+v)}$ was maximal in its (k+n-1)-th component) would have terminated the algorithm prior to iteration $i_3$ (where Step 4 of Section 4.3 is processed).

But, since the algorithm was not so terminated, we have to step back and state that $z_v^* \neq 0$. Thus, this subcase reduces to an instance of Subcase 1.1.

Case 2: $v = w$

Since the index set L contains at least two elements, there is another index, say $s$, in L such that $M^{(s)}$ has a positive component in row k+n+1. We can multiply $M^{(s)}$ by a positive number greater than 1 in order to make $M^{(s)}$ maximal in row k+n+1. Since $z_v^* + z_s^* > 0$, our proof of Case 1 applies in this case as well. Thus, this case reduces to an instance of Case 1.

The following table is for the $[M\ q]$ instance in iteration $i_3$, and it is intended to aid intuition regarding Claim 7.1; it shows $M^{(k+n+v)}$ as the $v$-th unit vector in iteration $i_3$.

| | $m_{1,w}$ | | $m_{1,v}$ | | | 0 | | $q_1$ | |
|---|---|---|---|---|---|---|---|---|---|
| | ⋮ | | ⋮ | | | ⋮ | | | |
| | | | | | | 0 | | | |
| | $m_{v,w}$ | | $m_{v,v}$ | | | 1 | | ⋮ | ← v |
| | | | | | | 0 | | | |
| | $m_{k+n,w}$ | | $m_{k+n,v}$ | | | ⋮ | | $q_{k+n}$ | |
| | +ve | | +ve & max | | | 0 | | $q_{k+n+1}$ +ve | |
| | | | ↑ v | | | ↑ k+n+v | | | |

END OF PROOF OF CLAIM 7.1

## 7.1.2 When a MinorP instance reverses a previous MajorP selection

We will state Claim 7.2 which is an analog of Claim 7.1. First, we present a two-part informal "background note".

.

Part 1 of background note - if (Eq) has a solution, then at least one unit-vector column of $M$ (in current $[M\,q]$) is not included in a solution basis matrix.

.

Details of this part of the note are as follows. Recall that $q_{k+n} = 0$, with some $q_j < 0$. We continue to assume that the columns of $[M\ q]$ are enumerated such that columns k+n+1,...,2(k+n) constitute the matrix U introduced earlier. If (Eq) has a solution, then one of the unit vectors (comprising U) corresponds to a $q_j < 0$, and it is not utilized in a solution of (Eq). A demonstration of that last statement is as follows.

If possible let $\hat{z}$ be a solution of (Eq) with $\hat{z}_{k+n+j} > 0$ for every $j$ having $q_j < 0$. Define a 2(k+n)-vector $\widetilde{z}$ by $\widetilde{z} \equiv (0,\ldots,0,q_1,\ldots,q_{k+n})^T$; note that some of the components of $\widetilde{z}$ are negative by virtue of some $q_j < 0$.

Clearly $M\widetilde{z} = q$. Let $\alpha$ be the smallest positive number such that $\alpha\hat{z} + \widetilde{z} \geq 0$. The vector $\alpha\hat{z} + \widetilde{z}$ has a 0 in at least one component, say the (k+n+t) component, with $q_t < 0$. Let $\omega$ denote the convex combination $(\alpha\hat{z} + \widetilde{z})/(\alpha + 1)$ of $\hat{z}$ and $\widetilde{z}$. The vector $\omega$ is non-negative, with 0 in at least one component corresponding to a $q_j < 0$. The vector $\omega$ is a solution of (Eq), because $M\omega = q$ &

$q_{k+n+1} = 0$; an LP complementary slackness theorem ensures that $\omega_j \omega_{(k+n+j)} = 0$, for $j = 1, \ldots, k+n$.

Thus, without any loss of generality, we can, in Claim 7.2 and its proof below, make the assumption that at least one of the unit-vector-columns corresponding to a $q_j < 0$ in $[M\, q]$ will not be utilized in the basis matrix of a solution of (Eq), if (Eq) has a solution.

To illustrate this background note, consider the following $[M\, q]$ instance that arises as the first $[M\, q]$ instance in solving the LP problem:

.

$$\begin{aligned} &\text{maximize} \quad 2x_1 - x_2 \\ &\text{subject to:} \quad x_1 + x_2 \leq 10 \\ &\qquad -x_1 \leq -1 \\ &\qquad x_1 \geq 0;\ x_2 \geq 0 \end{aligned}$$

.

that is,

$[M\, q] =$

| -10 | 1 | 3 | 0 | 1 | 0 | 0 | 0 | 10 | . |
|---|---|---|---|---|---|---|---|---|---|
| -10 | 1 | 1 | -1 | 0 | 1 | 0 | 0 | -1 | ←row 2 |
| -11 | 2 | 2 | -1 | 0 | 0 | 1 | 0 | -2 | ←row 3 |
| -11 | 1 | 2 | -1 | 0 | 0 | 0 | 1 | 1 | |
| -10 | 1 | 2 | -1 | 0 | 0 | 0 | 0 | 0 | |
| | | | | | ↑ 6 | ↑ 7 | | | |

In that $[M\, q]$, k=n=2, $q_2 < 0$, $q_3 < 0$, and $q_5 = 0$. It turns out that, after we perform two iterations of our algorithm on this example $[M\, q]$, we find $(2,0,10,0,0,9,0,3)^T$ to be a solution of this example's (Eq), thereby showing that $z_7 = 0$ as predicted in this background note.

.

`Part 2 of background note - if (Eq) has a solution, then at least two columns of` $M$ `(in current` $[M\, q]$`) are included in a solution basis matrix.`

.

Details of this part of the note are as follows. As $q_{k+n+1} = 0$ in current $[M\, q]$, if one column of $M$ having positive (k+n+1)-th component is included in a solution basis matrix, then another column of $M$ having negative (k+n+1)-th component must be included in that solution basis matrix as well. That is needed in order to have a balance in the (k+n+1)-th row of $[M\, q]$.

Thus, at least two columns of $M$ must be included in a solution basis matrix in this instance of $[M\, q]$. It is this observation that explains the difference between Step 4 of MinorP instance and Step 4 of MajorP instance descriptions in Section 4.3.

.

CLAIM 7.2*: Consider the algorithm processing Step 4 of "MinorP pivoting instance" in Section 4.3. If (Eq) has a solution, then there is a column of* $[M\, q]$*, denoted here as* $M^{(v)}$*, whose index v is contained in list L, such that either v and its complement are put into* $\Pi$ *or, else, the complementary GJ pivoting in* $M^{(v)}$ *terminates the algorithm with a solution of (Eq), as specified in Corollary 1.*

.

PROOF OF CLAIM 7.2

.

Assume that (Eq) has a solution, and recall details of Task 2 and Task 3 of Step 4 of "MinorP pivoting instance" in Section 4.3.

Utilizing "Part 1 of backgroung note", let $\widehat{j}$ be defined by column $\widehat{j}$ being a column of $M$ (in current $[M\ q]$) included in the basis matrix of a solution $z^*$ of (Eq), with $z^*_{k+n+\widehat{j}} = 0$.

The remainder of this proof consists of demonstrating that if Task 2 does not yield a solution of (Eq), then Task 3 will put $\widehat{j}$ and its complement into $\Pi$. We will demonstrate the contrapositive of that statement – if Task 3 would not put $\widehat{j}$ and its complement into $\Pi$, then Task 2 must have already yielded a solution of (Eq).

Recall that each iteration of the algorithm consists of a MinorP pivoting followed by a MajorP pivoting. Agreeing with the notation of Section 4.3's Step 4, let column $w$ be the column (of $M$) in which a MajorP pivoting would be done immediately after a MinorP pivoting has been done in column $\widehat{j}$.

Utilizing "Part 2 of background note", the (k+n+1)-th component of column $w$ is negative, whereas the (k+n+1)-th component of column $\widehat{j}$ is positive.

Next, as in Task 3 of Section 4.3's Step 4, consider reversing the order of GJ pivoting mentioned above, such that a complementary GJ pivoting is done in column $w$, followed by a complementary GJ pivoting being done in column $\widehat{j}$. Immediately after the complementary GJ pivoting in column $w$, resultant column $\widehat{j}$ has positive entry in its (k+n+1)-th component, along with at least one other column of $M$, say column $s$; that is on account of Task 3 not putting $\widehat{j}$ (and its complement) into $\Pi$.

Accordingly, we have $z^*_{\widehat{j}} + z^*_s > 0$; we also have $z^*_{k+n+\widehat{j}} = 0$ from "Part 1 of background note". Thus, the hypotheses of Lemma 6.1 and Corollary 1 are satisfied by column $\widehat{j}$ and column $s$. Therefore, we may next adopt our proof of Claim 7.1 in order to complete this proof of Claim 7.2. We therefore set $v \leftarrow \widehat{j}$ to have the conclusion of Claim 7.2.

.

END OF PROOF OF CLAIM 7.2.

.

## 7.2 Computational complexity lemma

We state and prove here a lemma, Lemma 7.1, about the proposed algorithm's computational complexity. Our proof of Lemma 7.1 is based on a tracking of MajorP column selections among $2(k + n - u)$ columns of M, where $u$ is the total number of times that the index set $\Pi$ is updated in the 4-step procedures described in Section 4.3.

.

LEMMA 7.1 *In at most* $2(k + n)$ *iterations, either the algorithm obtains a solution of (Eq) or, otherwise, it indicates that (Eq) has no solutions.*

.

PROOF OF LEMMA 7.1

.

First, note that if a run of the algorithm features a total of $u$ instances of L being a singleton, then the remainder of that run of the algorithm, outside of the $2u$ indices that comprise $\Pi$, has only $k + n - u$ complementary pairs of columns of $M$ left to be examined, because $M$ consists of $2(k + n)$ columns.

On account of an LP complementary slackness theorem, it is not necessary for the run of the algorithm to re-examine the $2u$ column indices that comprise $\Pi$, because, by the nature of complementary GJ pivoting, each column of $M$ is only competing (in the algorithm) against its own complement column.

There can be at most two MajorP column selections for each one of the $u$ column selection instances that comprise $\Pi$. Accordingly, $\Pi$ contributes *at most* $2u$ *iterations* to the algorithm's total number of iterations.

To count the remainder of the algorithm's iterations, there are two cases to consider, contingent upon whether some previous MajorP column selections are reversed/cancelled later in the algorithm – (i) no MajorP column selection gets reversed by a subsequent MajorP instance or MinorP instance, and (ii) a MajorP instance or a MinorP instance reverses a previous MajorP column selection, through Step 4 of the procedures described in Section 4.3.

Case (i): The case wherein no MajorP column selection gets reversed by a subsequent MajorP instance MinorP instance.

In this case, either a correct solution basis matrix is obtained before $k + n - u$ MajorP column selections have been made (in Step 3 of Section 4.3) or, otherwise, the algorithm indicates that (Eq) has no solutions by the instant that $k + n - u$ selections have been made (by MajorP instances in Step 3 of Section 4.3). This is true because $M$ has only $2(k + n - u)$ columns left to be examined, and each column selection by a MajorP instance implicitly reduces the number of "remaining" columns of $M$ by 2. Accordingly, at most $k + n - u$ iterations are added in this case, Case (i).

Case (ii): The case wherein a MajorP instance or a MinorP instance reverses a previous MajorP column selection through Step 4 of Section 4.3

That selection reversal must happen before $k + n - u$ number of MajorP selections have been made, as indicated in Case (i) above. Claims 7.1 and 7.2 above show that, if (Eq) has a solution, then such a MajorP selection reversal results in the algorithm being terminated, at that instant, with a solution of (Eq). The contrapositive of that statement is that if the index set L in Step 4 of Section 4.3 is exhausted without yielding a solution of (Eq), then it is because (Eq) has no solutions. Accordingly, at most $k + n - u$ iterations are added in this case, Case (ii).

Thus, as Case (i) and Case (ii) are mutually exclusive, the total number of iterations featured by each run of the algorithm is $2u + (k + n - u)$; that is, at most $2(k + n)$ iterations.

.

END OF PROOF OF LEMMA 7.1

# 8. Some directions for further work

Arising from this article, the following are some investigation topics that may be of interest to many members of the operations research or computing communities.

First, some work on the algorithm's numerical characteristics, especially its handling of ill-conditioned problem data, wherein extremely large numbers and extremely small numbers are mixed together, should ordinarily be of interest, in order to possibly make the algorithm useful in solving some classes of real-world LP problems.

Secondly, there may be useful connections between the algorithm and variants of the simplex algorithm, especially the primal-dual LP procedures, that may be utilized to enhance or explain the practical efficacy of some existing LP computer packages.

Finally, one may want to investigate how this algorithm can utilize what we call "MinorP column selection advice" along the way, without adversely affecting its computational complexity statement. The prospect for that is good because there is some flexibility for MinorP instances to select columns without reversing/cancelling previous MajorP column selections, as already indicated in Section 5 of this article. This could speed up this algorithm, and make it a useful sub-routine for solving classes of mixed-integer programming problems as well.

.

# 9. Appendix

.

On MinorP being well-defined

.

In Section 4.3, it is clear that each MajorP instance is well-defined, just on account of having $q_{k+n+1} > 0$. But it is not so clear that each MinorP instance is well-defined as well. With $q_{k+n+1} = 0$ and $q_i < 0$ some $i$, it is not clear that a column index $j$ having $m_{k+n+1,j} > 0$ exists.

We demonstrate through Claim 4 below that each MinorP instance is indeed well-defined, because there exists a column index $j$ having $m_{k+n+1,j} > 0$ so long as there is a row index $i$ with $q_i < 0$. Claim 4 and its proof accomplish that by demonstrating that, for $i = 1,\ldots,k+n$, either $m_{k+n+1,i}/q_i < 0$, or $m_{k+n+1,k+n+i}/q_i < 0$, for $q_i \neq 0$. As an implication of that statement, if $q_i < 0$, then either $m_{k+n+1,i} > 0$ or, else, $m_{k+n+1,k+n+i} > 0$.

.

Claim 4: *In Section 4.3, under MinorP pivoting instance description (that is, the $[M\ q]$ instance having $q_{k+n+1} = 0$), the ratio $m_{k+n+1,i}/q_i$ (or, else, $m_{k+n+1,k+n+i}/q_i$, when $m_{k+n+1,i} = 0$) has the same value for all row indices $i = 1,..,k+n$ having $q_i \neq 0$ in that $[M\ q]$ instance.*

.

Before we sketch a proof of Claim 4, we first give a numerical illustration.

.

*A numerical illustration:* Consider the following LP problem

$$
\begin{aligned}
\text{maximize:}\quad & 2x_1+x_2 \\
\text{such that:}\quad & x_1+x_2 \leq 5 \\
& x_1 \leq 2 \\
& x_1 \geq 0,\ x_2 \geq 0
\end{aligned}
$$

The first $[M\ q]$ instance, at Initialization, along with a last-column-last-row correspondence enumeration, is

.

| | | | | | | | | | |
|---|---|---|---|---|---|---|---|---|---|
| -5 | -2 | 3 | 2 | 1 | 0 | 0 | 0 | 5 | ←1 |
| -5 | -2 | 3 | 1 | 0 | 1 | 0 | 0 | 2 | ←2 |
| -6 | -3 | 2 | 1 | 0 | 0 | 1 | 0 | -2 | ←3 |
| -6 | -2 | 2 | 1 | 0 | 0 | 0 | 1 | -1 | ←4 |
| -5 | -2 | 2 | 1 | 0 | 0 | 0 | 0 | 0 | |
| ↑ | ↑ | ↑ | ↑ | | | | | | |
| 1 | 2 | 3 | 4 | | | | | | |

.

and the $m_{k+n+1,i}/q_i$ ratio is –1, for $i = 1,2,3,4$. After two complementary GJ pivotings, namely, MinorP at (4,4), followed by MajorP at (1,1), the $[M\ q]$, along with a last-column-last-row correspondence enumeration, is

.

| 1 | 0.29 | -0.14 | 0 | 0.14 | 0 | 0 | -0.29 | 1 | ←1 |
|---|---|---|---|---|---|---|---|---|---|
| 0 | -0.29 | 1.14 | 0 | -0.14 | 1 | 0 | -0.71 | 2 | ←2 |
| 0 | -1 | 0 | 0 | 0 | 0 | 1 | -1 | -1 | ←3 |
| 0 | -0.29 | 1.14 | 1 | 0.86 | 0 | 0 | -0.71 | 5 | ←4 |
| 0 | -0.29 | 0.14 | 0 | -0.14 | 0 | 0 | -0.71 | 0 | |
| | ↑ | ↑ | | ↑ | | | ↑ | | |
| | 2 | 3 | | 1 | | | 4 | | |

.

and the $m_{k+n+1,i}$ $/q_i$ (or else $m_{k+n+1,k+n+i}/q_i$) ratio here is -0.14. After two more complementary GJ pivotings, namely, MinorP at (3,3) after a copy of row 5 has been added to row 3 to enable the complementary pivoting in column 3, followed by MajorP at (2,2), the $[M\ q]$, along with a last-column-last-row correspondence enumeration, is

.

| 1 | 0 | 0 | 0 | 0.1 | 0.1 | 0.2 | -0.7 | 1 | ←1 |
|---|---|---|---|---|---|---|---|---|---|
| 0 | 1 | 0 | 0 | 0.1 | 0.1 | -0.8 | 1.3 | 1 | ←2 |
| 0 | 0 | 1 | 0 | -0.1 | 0.9 | -0.2 | -0.3 | 2 | ←3 |
| 0 | 0 | 0 | 1 | 1 | -1 | 0 | 0 | 3 | ←4 |
| 0 | 0 | 0 | 0 | -1 | -1 | -2 | -3 | 0 | |
| | | | | ↑ | ↑ | ↑ | ↑ | | |
| | | | | 1 | 2 | 3 | 4 | | |

.

and the $m_{k+n+1,k+n+i}/q_i$ ratio here is -1.

.

OUTLINE OF A PROOF OF CLAIM 4. We will state here an outline of our proof. A detailed proof is quite long, as one can see in our article entitled "On pairs of complementary Gauss-Jordan pivotings transforming skew-symmetric matrices" posted at arxiv.org as arXiv:2410.19350.

At initialization, the matrix $[M\ q]$ is of the form

$[M\ q] =$

| 0 | $m_{12}$ | ⋯ | $m_{1,k+n}$ | 1 | 0 | ⋯ | 0 | $q_1$ |
|---|---|---|---|---|---|---|---|---|
| $-m_{12}$ | 0 | ⋯ | $m_{2,k+n}$ | 0 | 1 | | ⋮ | $q_2$ |
| ⋮ | ⋮ | ⋱ | ⋮ | ⋮ | | ⋱ | 0 | ⋮ |
| $-m_{1,k+n}$ | $-m_{1,k+n}$ | ⋯ | 0 | 0 | ⋯ | 0 | 1 | $q_{k+n}$ |
| $-q_1$ | $-q_2$ | ⋯ | $-q_{k+n}$ | 0 | 0 | 0 | 0 | 0 |

with the submatrix

$$S = \begin{array}{|c|c|c|c|c|}
\hline
0 & m_{12} & \cdots & m_{1,k+n} & q_1 \\
\hline
-m_{12} & 0 & \cdots & m_{2,k+n} & q_2 \\
\hline
\vdots & \vdots & \ddots & \vdots & \vdots \\
\hline
-m_{1,k+n} & -m_{1,k+n} & \cdots & 0 & q_{k+n} \\
\hline
-q_1 & -q_2 & \cdots & -q_{k+n} & 0 \\
\hline
\end{array}$$

being skew-symmetric. Our proof, fully described in arXiv:2410.19350, is a mathematical induction proof showing how a sequence of "pairs of complementary GJ pivotings", corresponding to MinorP and MajorP, transform S into a sequence of skew-symmetric matrices.

.

END OF CLAIM 4 PROOF OUTLINE

.

.

# 10. References

1. Bland, R. G., Goldfarb, D., Todd, M. J., (1981) "The ellipsoid method: a survey," *Oper. Res*. 29, 1039-1091

2. Chubanov, S. (2011) "A strongly polynomial algorithm for linear systems having binary solutions" *Math Programming*, p.1-38

3. Cottle, R. W., J.-S. Pang, R. E. Stone (2013) *The Linear Complementarity Problem,* Academic Press

4. den Hertog, D., C. Roos, T. Terlaky (1993) "The linear complementarity problem, sufficient matrices, and the criss-cross method" *Linear Algebra and its Applications*, 187, 1-14 , submitted by R. Cottle

5. Dantzig, G. B., Thapa, M. N. (1997) *Linear Programing, 1: Introduction*, Springer, New York

6. Dantzig, G. B., Thapa, M. N. (2003) *Linear Programing, 2: Theory and Extensions*, Springer, New York

7. Goldfarb, D., Todd, M. J., "Linear programming," in: G. L. Nemhauser, A. H. G. Rinnooy Kan and M. J. Todd, eds., (1989) *Optimization, Handbooks in Operations Research and Management Science*, Vol. I, North-Holland, Amsterdam, 73-170

8. Kantorovich, L. V. (1939) "Mathematical methods in the organization and planning of production," [English translation: *Management Science* 6, (1960) 366-422]

9. Karmarkar, N. (1984) "A new polynomial-time algorithm for linear programming," *Combinatorica* 4, 373-395

10. Khachian, L. G. (1979) "A polynomial algorithm for linear programming, "*Doklady Akademiia* Nauk USSR 244, 1093-1096 [English translation: Soviet Mathematics Doklady 20, 191-194

11. Kojima, M., S. Misuno, A. Yosishe (1999) "A polynomial-time algorithm for a class of linear complementarity problems," *SIAM Journal on Optimization*, vol 10, #1

12. Megiddo, N. (1984) "Linear programming in linear time when the dimension is fixed," *Journal of the Association of Computing Machinery* 31, 114-127

13. Nering, E. D., A. W. Tucker (1993) *Linear Programs and Related Problems*, Academic Press

14. Saigal, R. (1995) *Linear Programming A Modern Integrated Analysis*, Kluwer Academic Publishers

15. Tardos, E. (1986) "A strongly polynomial algorithm to solve combinatorial linear programs," *Operations Research* 35, 250-256

16. Todd, M. J. (1988) "Polynomial algorithms for linear programming," in *Advances in Optimization and Control*, H. A. Eiselt and G. Pederzoli, eds., Springer-Verlag, Berlin,

17. Wheaton, I., S. Awoniyi (2017) "A new iterative method for solving non-square systems of linear equations," *Journal of Computational and Applied Math*. 322, 1-6

18. arXiv:2503.12041, Awoniyi, S. A., (March 2025) "A strongly polynomial-time algorithm for the general linear programming problem"

19. arXiv:2410.19350, Awoniyi, S. A., (October 2024) "On pairs of complementary GJ pivoting transforming skew-symmetric matrices"